\documentclass[reqno,11pt]{amsart} 
\usepackage[top=1.5in,right=1.125in,left=1.125in,bottom=1.5in]{geometry}
\usepackage{amssymb}
\usepackage{color}
\usepackage{tikz}
\usetikzlibrary{positioning,decorations.pathmorphing}
\usetikzlibrary{positioning,decorations.pathmorphing}
\definecolor{red}{rgb}{0.7,0,0}
\definecolor{grey}{RGB}{112,112,112}
\definecolor{blue}{RGB}{034,113,179}
\usepackage[colorlinks=true,citecolor=blue,linkcolor=red,urlcolor=grey]{hyperref}

\theoremstyle{remark}

\newcounter{mnotecount}[section]

\renewcommand{\themnotecount}{\thesection.\arabic{mnotecount}}

\newcommand{\mnote}[1]
{\protect{\stepcounter{mnotecount}}$^{\mbox{\footnotesize
$
\bullet$\themnotecount}}$ \marginpar{
\raggedright\tiny\em
$\!\!\!\!\!\!\,\bullet$\themnotecount: #1} }

\newcommand{\CP}{\mathbb{CP}}

\newcommand{\C}{\mathbb{C}}

\newcommand{\R}{\mathbb{R}}

\newcommand{\pr}[2]{\langle#1, #2\rangle}

\def\be{\begin{equation}}

\def\ee{\end{equation}}

\def\bea{\begin{eqnarray}}
\def\eea{\end{eqnarray}}

\numberwithin{equation}{section}
\begin{document} \date{}
\title{Twistor theory of Loxodromes}
\author{Maciej Dunajski}
\address{Department of Applied Mathematics and Theoretical Physics\\ 
University of Cambridge\\ Wilberforce Road, Cambridge CB3 0WA, UK.}
\email{m.dunajski@damtp.cam.ac.uk}
\author{Wojciech Kry\'nski}
\address{
Institute of Mathematics, Polish Academy of Sciences, ul. Sniadeckich 8, 00-656 Warszawa, ´
Poland.}
\email{krynski@impan.pl}
\maketitle
\begin{center}
{\em Dedicated to Roger Penrose
on the occasion of his 95th birthday.}
\end{center}
\begin{abstract}
The twistor correspondence provides a duality between curves in $S^4$ and ruled surfaces in $\CP^3$ in which the conformal properties
of the former are reflected in the projective properties of the latter.
We use this to characterise a $14$--dimensional family of homogeneous ruled surfaces in $\CP^3$ which 
correspond to loxodromes in $S^4$.
\end{abstract}
\section{Introduction}
A loxodrome, or rhumb line, is a curve on a sphere that crosses all meridians at a constant
angle $\psi$. Historically, marine navigation along a rhumb line was easier than continually
adjusting the compass bearing, as would be needed along a great--circle geodesic. The price paid for
this  is modest: a loxodrome winds infinitely often around each pole and yet has finite
length, exceeding that of the meridian arc between the same two latitudes by a factor of
$\sec\psi$. The Mercator projection from a sphere to a plane yields a 
map\footnote{The loxodrome was identified as a curve
distinct from the great circle by Pedro Nunes in 1537. 
 Mercator's world map of 1569
exploited the fact that a suitable conformal projection straightens it.  Mercator has  never explained
his construction, and the mathematics was supplied by Edward Wright in 1599. See \cite{synder}.} 
on which loxodromes are straight lines. This projection is not an isometry, but
it is conformal, as it preserves angles. Another conformal map from spheres to planes is the stereographic
projection. It takes loxodromes to logarithmic spirals centred at the image of a pole. The two poles are the limit
points of a loxodrome.
\begin{center}
    \includegraphics[width=9cm,keepaspectratio]{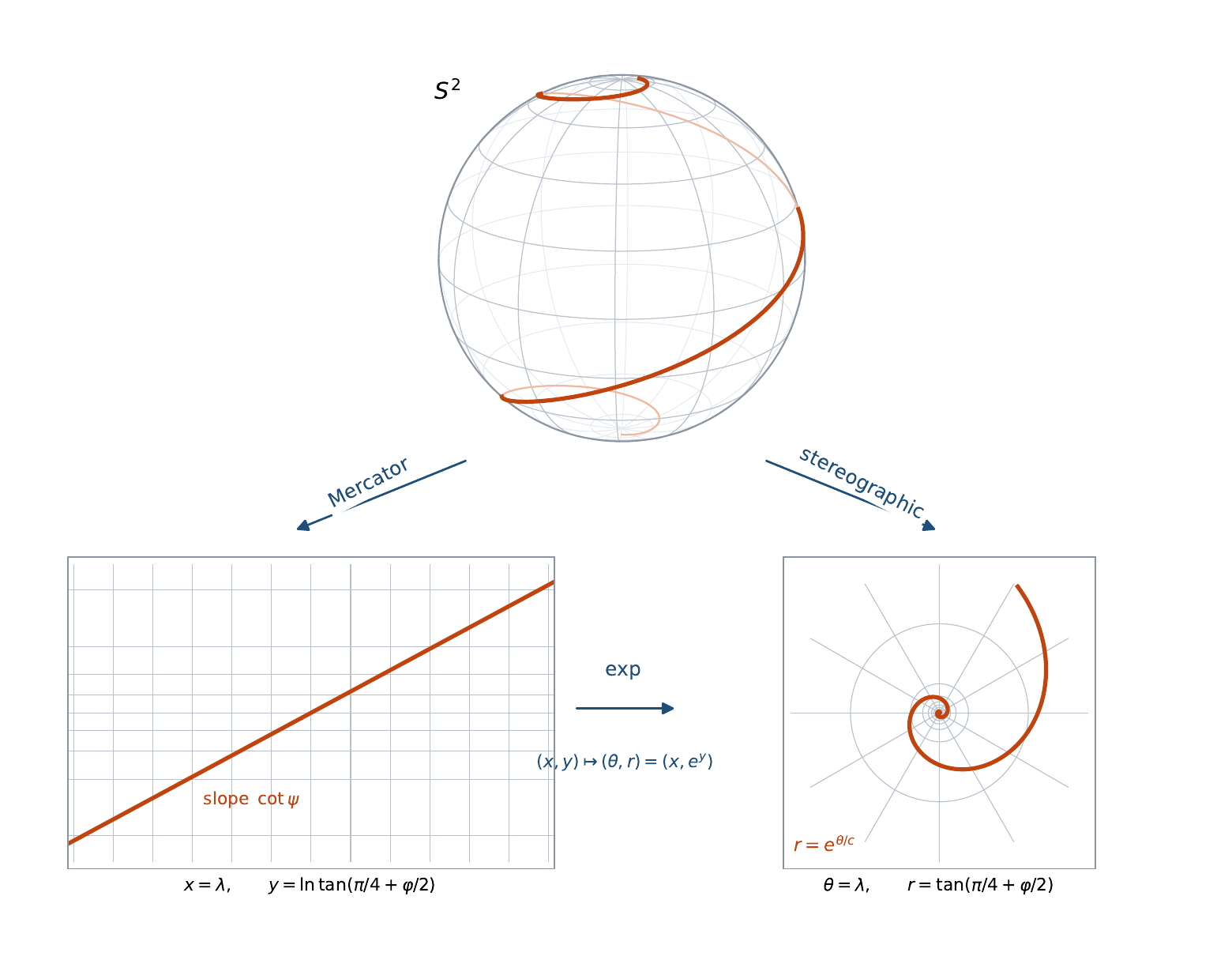}
\end{center}
In this  paper we shall use the twistor theory of Roger Penrose \cite{Pe78, AHS} to explore the conformal differential geometry of loxodromes on spheres,
regarded as integral curves of a system of 4th order ordinary differential equations.  Twistor theory aspires to unify classical and quantum physics, and employing it to
address a problem of 16th--century marine navigation sounds like using a sledgehammer to crack
a nut. But the sledgehammer is an attractive one and we shall show in \S\ref{mainresult} that it reveals
 properties which may not have been noticed before: In twistor theory the points
on the four--sphere correspond to certain projective lines in $\CP^3$. We shall employ this 
to characterise a $14$--dimensional family of homogeneous ruled surfaces in $\CP^3$ which 
correspond to loxodromes in a way which swaps the conformal invariance with the projective invariance.

\subsection*{Acknowledgements} 
This work has began when the authors attended 
the {\em Conformal Geodesics}  programme held at the Banff International Research Station in August 2024.
We are grateful to the participants of the programme (in particular to Mike Eastwood) for discussions. 
Anthropic's  Opus 5 was used to generate Figures.
This research
was supported by the Simons Foundation grant SFI-MPS-T-Institutes-00010825, and by the
State Treasury funds under the project Organization of the Simons Semesters at the Banach Center
- New Energies in 2026-2028 (MNiSW/2025/DAP/491).

\section{Twistor correspondence}
Twistor theory  gives
a correspondence between the lines in the twistor space $\CP^3$ and the points in the four--dimensional Grassmannian $\mbox{Gr}_2(\C^4)$ which in the original construction of Penrose is identified with the complexified and compactified Minkowski space. The conformal geometry of $\mbox{Gr}_2(\C^4)$ is encoded in the projective geometry of $\CP^3$: two points in the Grassmannian are separated 
by a null curve iff the corresponding lines intersect
in $\CP^3$. The real version of this construction
due to Atiyah, Hitchin and Singer \cite{AHS} singles out {\em real} lines in $\CP^3$ which correspond to points in the four--dimensional sphere $S^4$. The resulting conformal structure, when restricted
to the real slice, containts the round metric on $S^4$. Moreover the
$PSL(4, \C)$ projective transformations of $\CP^3$  which map lines to lines correspond to conformal transformations of the Grassmannian or the four--sphere.

It follows that surfaces
in $\CP^3$ ruled by lines are in one--to--one correspondence with curves
in the Grassmannian, and so
the classification of homogeneous curves in conformal geometry  corresponds to a classification of certain ruled surfaces in projective three--space.  The conformal circles correspond to quadrics  which are doubly ruled \cite{BE} . This leads to interesting dualities. For example the existence of a  
unique circle through three points  is dual to the statement that  there is a unique quadric containing three skew lines. 

In \S\ref{mainresult} we shall identify a class of homogeneous ruled surfaces in $\CP^3$
invariant under the anti--holomorphic involution
which, via the twistor transform, correspond to  the conformal Loxodromes \cite{DK, E23, lenka}.
\subsection{The incidence relation}
The twistor space is defined to be the complex projective 3-space ${\bf T}=\CP^3$. Let
$[X, Y, W, Z]\sim [\rho X, \rho Y, \rho W, \rho Z]$ where $\rho\in \C^{\star}$ be homogeneous
complex coordinates on this space. The link between ${\bf T}$ and $\mbox{Gr}_2(\C^4)$ can be concretely expressed by removing a projective line $[X, Y, 0, 0]$ from ${\bf T}$ and considering
the complexified Minkowski space $\C^4$ inside the Grassmannian with coordinates $(x^1, x^2, x^3, x^4)$.
The incidence relation between the two spaces is then \cite{Pe78, ADM}
\be
\label{twistorequation}
X=W(x_1+ix_2)+(x_3+ix_4)Z, \quad Y=-W(x_3-ix_4)+Z(x_1-ix_2).
\ee
Fixing a point in $\C^4$ gives a projective line $\CP^1$ in ${\bf T}$ and conversely fixing a point
in ${\bf T}$ gives a totally self--dual plane (the $\alpha$--plane in the terminology of \cite{Pe78}) in $\C^4$.
\begin{eqnarray*}
\mbox{Complexified space-time $\C^4$}\quad&\longleftrightarrow&\mbox{Twistor space ${\bf{T}}$}.\\
\mbox{Point $p$}\quad&\longleftrightarrow&\mbox{Complex line}.\,\, \
L_p=\CP^1\\
\mbox{Null self-dual (=$\alpha$) two-plane }\quad&\longleftrightarrow&\mbox{Point}.\\
p_1, p_2 \;\mbox{null separated}  \quad&\longleftrightarrow&
L_1, L_2 \; \mbox{intersect at one point.}
\end{eqnarray*}
The resulting conformal structure is deduced from (\ref{twistorequation}) and gives a cone
${(x^1)}^2+\dots +{(x^4)}^2=0$.
The conformal transformations of the Grassmannian which preserve this flat conformal structure correspond to
projective transformations of ${\bf{T}}$: the conformal group of complexified space--time is $SL(4,\C)$.
\subsection{Reality conditions} 
The twistor lines which are preserved by an anti--holomorphic involution $\tau:\CP^3\rightarrow \CP^3$
\be
\label{holinvolution}
\tau [X, Y, W, Z]=[-\bar{Y}, \bar{X}, -\bar{Z}, \bar{W}]
\ee
form a four--dimensional real family $\R^4$, or (if a line has not been removed from ${\bf T})$ $S^4$.
This involution does not have any fixed points in ${\bf T}$, so for any $\zeta\in {\bf T}$ the two points
$(\zeta, \tau(\zeta))$ are connected by a unique real line. This leads to the Atiyah--Hitchin--Singer version of the twistor correspondence \cite{AHS} where ${\bf T}$ is regarded as a real six--dimensional manifold which fibers over $S^4$. This fibration is smooth, but not holomoprhic and a fiber over $p\in S^4$ is the twistor line 
$\CP^1$.
The projective transformations preserving the 
real lines reduce $SL(4, \C)$ to $SL(2,\mathbb{H})$,
which in turn is the  double cover of the conformal group of $S^4$. In fact, under the identification
$S^4=\mathbb{HP}^1$, the conformal symmetries of $S^4$ are represented by quaternionic M\"obius
transformations.
\section{Conformal Mercator equation}
In this section we use capital letters to denote vectors in $\R^4$, and $\pr{X}{Y}$ is the inner product of two vectors with respect to a flat metric. We also set $|X|^2=\pr{X}{X}$. If $X=X(t)$ then $U=\dot{X}$ and $A=\dot{U}$.

The conformally invariant Mercator equation \cite{DK, lenka} is a system of 4th order ODEs arising from a conformally invariant second order Lagrangian. To this end we shall restrict to the flat conformal structure on $\mbox{Gr}_2(\C^4)$ and
write the corresponding equation using the flat metric in this conformal class as 
\be
\label{4th1}
\frac{dC}{dt}=0 \quad\mbox{where}\quad
 C=\frac{1}{|U|^2}\Big(\dot{A}-\frac{|A|^2}{|U|^2}U-2\frac{\pr{A}{U}}{|U|^2}A
       +4\frac{\pr{A}{U}^2}{|U|^4}U
-2\frac{\pr{\dot{A}}{U}}{|U|^2}U \Big). 
\ee
A solution curve to \eqref{4th1} is determined by specifying initial conditions $X(0), U(0), A(0), \dot{A}(0)$. These initial conditions determine  $C$ in \eqref{4th1}, so one can instead use $X(0), U(0), A(0), C$ as initial condition as these determine $\dot{A}(0)$. 

All conformal circles are integral curves of \eqref{4th1} for special initial conditions. To see it set
\be
\label{C_formula}
C=\Big(\frac{1}{2}\frac{|A|^2}{|U|^4}-2\frac{\pr{A}{U}^2}{|U|^6}\Big)U+\frac{\pr{A}{U}}{|U|^4}A,
\ee
and substitute this into \eqref{4th1}. This yields the conformal geodesic equations for the flat metric (compare \cite{BE})
\be
\label{flateq}
\dot{A}-\frac{3\pr{A}{U}}{|U|^2}A+\frac{3|A|^2}{2|U|^2}U=0
\ee
with projectively parametrised circles
\be
\label{circles1}
t\rightarrow X(t)=X_0+\frac{tU_0+t^2 A_0}{1+t^2|A_0|^2},
\ee
as solution curves. Here
$U_0$ is a constant unit vector, and $\pr{U_0}{A_0}=0$.

For generic initial conditions the integral curves of \eqref{4th1} are not conformal geodesics.  For example, choosing arbitrary values of $X(0), U(0), A(0)$, and setting $C=0$ reduces \eqref{4th1} to
\be
\label{seq11}
\dot{A}-\frac{2\pr{A}{U}}{|U|^2}A+\frac{|A|^2}{|U|^2}U =0.
\ee
The general solution of this system
is
\be
\label{spirals}
t\rightarrow X(t)=e^t\cos{(ct)}\;P_0 +e^t\sin{(ct)}\;Q_0+R_0,
\ee
where $P_0, Q_0, R_0$ are constant vectors such that $\pr{P_0}{Q_0}=0$ and $|P_0|=|Q_0|$. The curves \eqref{spirals} are logarithmic spirals in the plane spanned by $(P_0, Q_0)$ which spiral towards $R_0$ as $t\rightarrow -\infty$. If $(r, \theta)$ are plane polar coordinates in the plane spanned by $P_0, Q_0$ and centered at $R_0$,  then the unparametrised form of the spirals\footnote{
It was an open problem \cite{Tod} to determine whether spirals can occur as conformal geodesics. This has now been settled affirmatively by Wojciech Kami{\'n}ski \cite{Kaminski}.
}
is $r=|P_0|e^{\theta/c}$ where $c$ is related to the constant angle $\psi$ from the Introduction by $c=\tan{\psi}$.

The conformal invariance of the fourth order system \eqref{4th1} ensures  that the inverse images of the logarithmic spirals under the stereographic projection from $S^{4}$ to $\R^4$ are solutions to \eqref{4th1} on the round sphere. These curves are the loxodromes. They cut all meridians at a fixed angle, and correspond to straight lines on the Mercator map which justifies our terminology.

 The general solution to the conformal Mercator equation \eqref{4th1} is given by the special conformal transformation of the logarithmic spiral \eqref{spirals}:
\be
\label{final_s}
t\rightarrow Y(t)=\frac{X(t)-|X(t)|^2B}{1-2\pr{X(t)}{B}+|B|^2 |X(t)|^2},
\ee
where $X(t)$  is given by \eqref{spirals}, and $B$ is a constant vector.

The spirals corresponding to \eqref{seq11} can be made to spiral to two arbitrary spiralling points which  can be constructed  by a conformal transformations. The action the conformal group on the solutions to the 4th order equation \eqref{4th1} has two orbits. The circles belong to one orbit, and the spirals to another containing $C=0$.
 
\subsection{Null geodesics}
Null geodesics form a subset of circles.  The corresponding surfaces are flags consisting of a point $\zeta\in{\bf T}$ and a plane $\Pi$ in $\CP^3$ such that $\zeta\in \Pi$. Any such flag is also a ruled surface: it is ruled by all lines in $\Pi$ which pass through $\zeta$.

Indeed, any null geodesic is an intersection of an 
$\alpha$ surface and a $\beta$ surface in $\mbox{Gr}_2(4)$. The $\alpha$ surface corresponds to a point in $\CP^3$, and the $\beta$ surface
to a point in ${\CP^3}^{*}$, or equivalently to a plane in $\CP^3$. The incidence between the two 
gives the flag. This flag depends
on $3$ (a point in $\CP^3$) $+$ $3$ (a point in 
${\CP^3}^*$) $-1$ (incidence) $=5$ parameters, which is the dimension of the space of null geodesics.

\section{Ruled homogeneous surfaces}
\label{mainresult}
Not all solutions of the Mercator equation belong to one orbit of the conformal group. While $P_0, Q_0, R_0, B$ can be adjusted by the conformal transformations, the parameter $c$ can not - it is conformally invariant and fixes the angle $\tan^{-1}(c)$ at which the loxodromes cut the meridians. Therefore the unparametrised integral curves of the Mercator equation form a 14--dimensional 
family. This agrees with the following count: For a 4th order equation we specify
$X, \dot{X}, \ddot{X}, \dddot{X}$ at $t=0$ which gives 16 parameters. The equation is invariant under affine reparametrisation of $t$ and subtracting $2$ gives $14=16-2$ (the projective freedom is enjoyed only by the sub-family of conformal circles).
The $SL(4, \C)$ orbits are 13--dimensional (or 9--dimensional if they are circles). Thus, in the twistor space, they correspond to homogeneous projective surfaces with two--dimensional symmetry 
($2=\mbox{dim}_{\C}(SL(4, \C))-13$).
In this section we shall find these surfaces, show that they are preserved by the  anti--holomorphic involution (\ref{holinvolution}), and identify them in the classification of \cite{Dillen, Doubrov}.

Consider the spiral (\ref{spirals})
\[
x=(x_1, x_2, x_3, x_4)=P_0 e^t\cos{ct}+Q_0 e^t\sin{ct}, \quad
|P_0|=|Q_0|, \quad \pr{P_0}{Q_0}=0
\]
and  substitute this into the twistor incidence relation (\ref{twistorequation}).
Choosing
$
P_0=(0, 0, 1, 0),  Q_0=(1, 0, 0, 0), R_0=(0, 0, 0, 0)
$
yields 
\[
X=We^{t}\sin{ct}+Z e^{t}\cos{ct}, \quad Y=-We^{t}\cos{ct}+Ze^{t}\sin{ct}.
\]
Eliminating $t$ between these two equations gives
\be
\label{surface1}
\ln\frac{X^2+Y^2}{W^2+Z^2}=\frac{2}{c}\arctan{\frac{WX+ZY}{ZX-WY}},
\ee
or, using the affine chart with $W=1$, 
\be
\label{surface2}
\ln\frac{X^2+Y^2}{1+Z^2}=-\frac{2}{c}\Big(\arctan(X/Y)+\arctan{Z}\Big)
\ee
which we recognize as a one--parameter family of homogeneous surfaces $|20|$ in \cite{Dillen} 
or  family $(3)$ in \cite{Doubrov}.
For each $c$ the surface (\ref{surface1}) admits a two-dimensional symmetry subgroup of
$SL(4, \C)$.
The fourteen--parameter family of ruled surfaces in the twistor space $\CP^3$ corresponds is a family
of $13$--dimensional orbits of $SL(4, \C)$ with $c$  as the additional parameter.
\begin{center}
    \includegraphics[width=12cm,keepaspectratio]{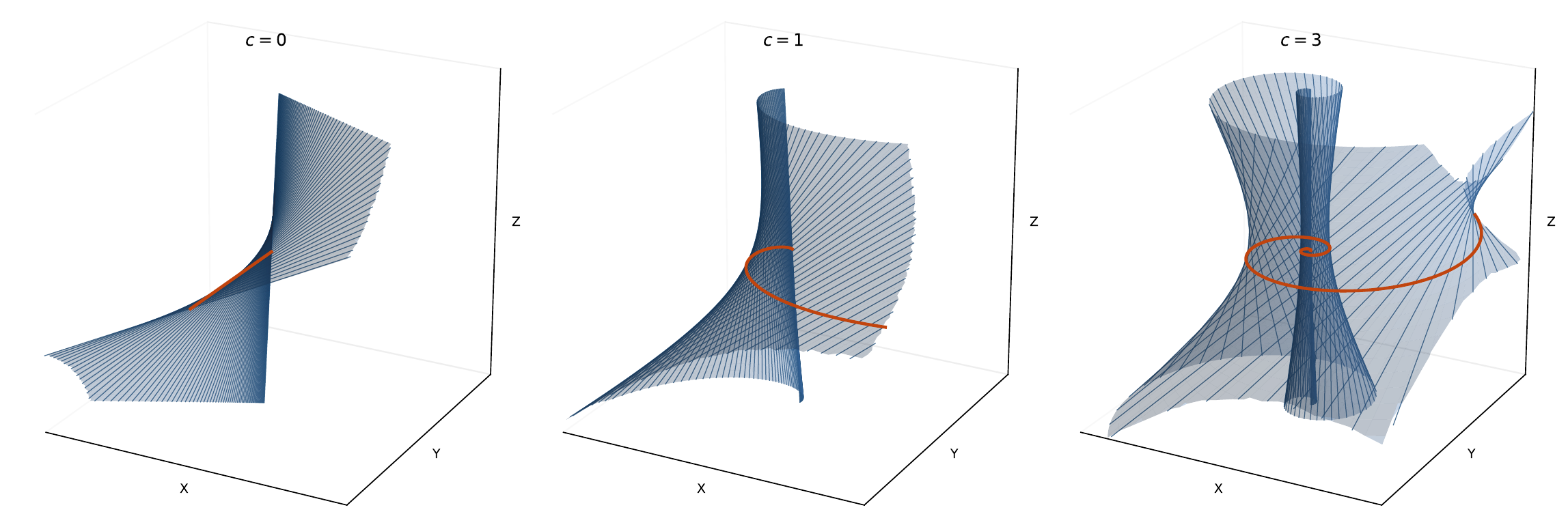}
\end{center}
\begin{center}
{\em  Selected branches of surfaces (\ref{surface2}) in the real affine chart $W=1$.}
\end{center}
\begin{center}
    \includegraphics[width=6cm,keepaspectratio]{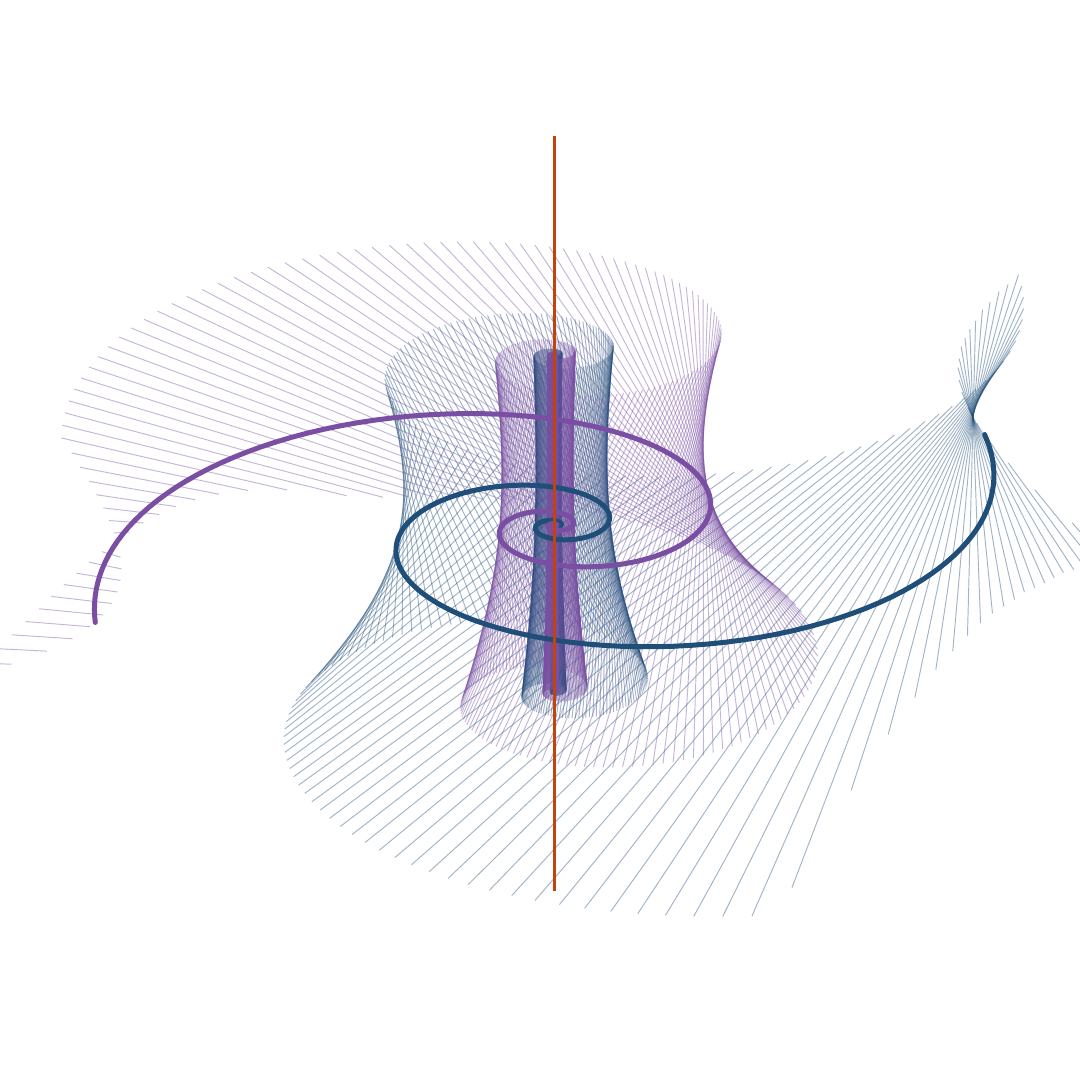}
\end{center}
\begin{center}
{\em  Two branches $k=0$ and $k=1$ of $\ln\frac{X^2+Y^2}{1+Z^2}=-\frac{2}{c}\Big(\arctan(X/Y)+\arctan{Z}+k\pi\Big)$. }
\end{center}
\subsection{Reality conditions}
The surface (\ref{surface1}) is invariant under the anti--holomorphic involution (\ref{holinvolution}).
The twistor lines invariant under this involution correspond to points in $\R^4$, or $S^4$.
\subsection{Conformal circles}
 In the special case
$c=0$ the spirals reduce to conformal geodesics (circles) and 
the surface (\ref{surface1}) reduces to a quadric
\[
WX+ZY=0.
\]
The $SL(4, \C)$ orbit of this quadric is $9$--dimensional, in agreement with the statement in \cite{BE}
that conformal circles on complexified $S^4$  correspond to doubly--ruled quartics in the twistor space.
\subsection{Ruled cubic surfaces}
The surface  (\ref{surface1}) can be represented as
\[
\Big(\frac{X+iY}{Z-iW}\Big)^{1-ic}=\Big(\frac{X-iY}{Z+iW}\Big)^{1+ic}.
\]
Setting $c=ip/q$ where $p, q$ are co-prime integers and applying $SL(4, \C)$ yields an projectively equivalent surface
\be
\label{surfacerat}
X^{q+p}W^{q-p}=Y^{q-p}Z^{q+p}
\ee
which is algebraic and of degree $2q/\mbox{gcd}(q+p, q-p)$. If $c=i/3$ then (\ref{surfacerat}) is a ruled cubic surface 
$X^2W=YZ^2$
with 2--dimensional symmetry (so it is different than the Cayley ruled cubic which has 3--dimensional symmetry). The rational form 
(\ref{surfacerat}) is not preserved by the involution (\ref{holinvolution}) so there are no corresponding 
real curves on $\R^4$ or $S^4$. There are such curves on Minkowski space forming a subclass of all parabolas.

\end{document}